\documentclass[12pt]{amsart}
\usepackage{amsfonts, amstext, amsmath, amsthm, amscd, amssymb}
\usepackage{color, graphicx, hyperref,  verbatim, subfigure, rotating, multirow, wasysym, float} 
\let\oldmarginpar\marginpar
\renewcommand\marginpar[1]{\oldmarginpar[\raggedleft\footnotesize #1]
{\raggedright\footnotesize #1}}
\newtheorem{defn}{Definition}

\newtheorem{thm}{Theorem}

\newtheorem{lem}{Lemma}

\title[On 3-braids and L-space knots.]{On 3-braids and L-space knots}
\author[C. Lee]{Christine Ruey Shan Lee} 
\address{Department of Mathematics and Statistics, University of South Alabama, Mobile, AL 36608} 
\email{christine.rs.lee@gmail.com}

\author[F. Vafaee]{Faramarz Vafaee} 
\address{Department of Mathematics, Duke University, Durham, NC 27708} 
\email{vafaee@math.duke.edu}
\begin{document}

\maketitle

\begin{abstract}
We classify closed 3-braids which are L-space knots. 
\end{abstract}

\section*{Introduction} A rational homology 3-sphere $Y$ is an L-space if $|H_1(Y; \mathbb{Z})|=\text{rank } \widehat{HF}(Y)$, where $\widehat{HF}$ denotes the `hat' version of Heegaard Floer homology, and the name stems from the fact that lens spaces are L-spaces. Besides lens spaces, examples of L-spaces include all connected sums of manifolds with elliptic geometry \cite{Ath}. 

A prominent source of L-spaces arises from surgeries on knots. Suppose that K is a knot
in $S^3$: if K admits a non-trivial surgery to an L-space, then K is an {\emph L-space knot}. Examples include torus knots and, more generally, Berge knots in $S^3$. Various properties of L-space knots have been studied in the previous years; the two particularly pertinent to our work are about the Alexander polynomial $\triangle_K(t)$ of an L-space knot $K$: 
\begin{itemize}
\item The absolute value of a nonzero coefficient of $\triangle_K(t)$ is 1. The set of nonzero coefficients alternates in sign \cite[Corollary~1.3]{Ath}. 
\item If $g$ is the maximum degree of $\triangle_K(t)$ in $t$, then the coefficients of the term $t^{g-1}$ is nonzero and therefore $\pm 1$ \cite{Hedden2014}. 
\end{itemize} 

The purpose of this manuscript is to study which 3-braids, that close to form a knot, admit L-space surgeries. We prove that:
\begin{thm}
Twisted $(3, q)$ torus knots are the only knots with 3-braid representatives that admit L-space surgeries.
\label{thm:2}
\end{thm}
Our proof of Theorem~\ref{thm:2} uses the constraints on the Alexander polynomial of an L-space knot that have previously been studied to give the classification of L-space knots among pretzel knots \cite{li-mo}. %(This was later generalized to We show 
We begin by computing certain coefficients of the Jones polynomials of closed 3-braids.
%In this paper, we classify the L-space knots among closed 3-braids by using the Jones polynomial \cite{farey}. 
The Alexander polynomial of a closed 3-braid may be written in terms of the Jones polynomial \cite{birman}. This allows us to compute certain coefficients of the Jones polynomials using \cite{farey} to rule out closed 3-braids whose Alexander polynomials violate the aforementioned conditions.  \\
       
\paragraph{\textbf{Acknowledgement.}} We would like to thank Efstratia Kalfagianni for referring us to the problem. We would also like to thank Cameron Gordon and Kenneth Baker for their conversations with the first author and their interest in this work. 

\section*{The Jones polynomial, 3-braids, and the Alexander polynomial}	

We will first derive an expression of the Alexander polynomial of a closed 3-braid in terms of the Jones polynomial. Let $B_n$ be the $n$-string braid group. The \emph{Burau representation} of $B_n$ is a map $\psi$ from $B_n$ to $n-1 \times n-1$ matrices with entries in $\mathbb{Z}[t, t^{-1}]$. 
\[ \psi: B_n \rightarrow GL(n-1, \mathbb{Z}[t, t^{-1}]).  \] 

\begin{figure}
\centering
\def \svgwidth{.4\columnwidth}
%% Creator: Inkscape 1.0beta1 (32d4812, 2019-09-19), www.inkscape.org
%% Accompanies image file 'generator.pdf' (pdf, eps, ps)
\begingroup%
  \makeatletter%
  \providecommand\color[2][]{%
    \errmessage{(Inkscape) Color is used for the text in Inkscape, but the package 'color.sty' is not loaded}%
    \renewcommand\color[2][]{}%
  }%
  \providecommand\transparent[1]{%
    \errmessage{(Inkscape) Transparency is used (non-zero) for the text in Inkscape, but the package 'transparent.sty' is not loaded}%
    \renewcommand\transparent[1]{}%
  }%
  \providecommand\rotatebox[2]{#2}%
  \newcommand*\fsize{\dimexpr\f@size pt\relax}%
  \newcommand*\lineheight[1]{\fontsize{\fsize}{#1\fsize}\selectfont}%
  \ifx\svgwidth\undefined%
    \setlength{\unitlength}{277.10134191bp}%
    \ifx\svgscale\undefined%
      \relax%
    \else%
      \setlength{\unitlength}{\unitlength * \real{\svgscale}}%
    \fi%
  \else%
    \setlength{\unitlength}{\svgwidth}%
  \fi%
  \global\let\svgwidth\undefined%
  \global\let\svgscale\undefined%
  \makeatother%
  \begin{picture}(1,0.49105004)%
    \lineheight{1}%
    \setlength\tabcolsep{0pt}%
    \put(0,0){\includegraphics[width=\unitlength,page=1]{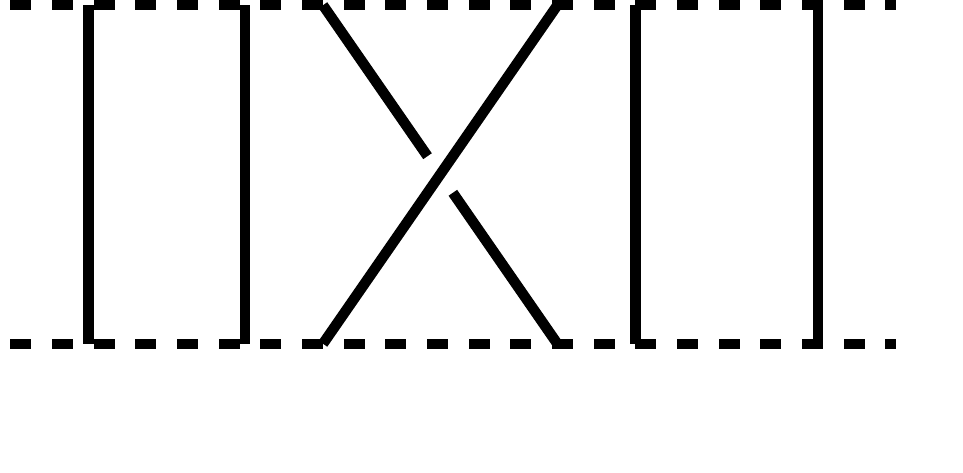}}%
    \put(0.05676449,0.01908352){\color[rgb]{0,0,0}\makebox(0,0)[lt]{\lineheight{1.25}\smash{\begin{tabular}[t]{l}$1$\end{tabular}}}}%
    \put(0.29780944,0.01908352){\color[rgb]{0,0,0}\makebox(0,0)[lt]{\lineheight{1.25}\smash{\begin{tabular}[t]{l}$i$\end{tabular}}}}%
    \put(0.82830117,0.01908352){\color[rgb]{0,0,0}\makebox(0,0)[lt]{\lineheight{1.25}\smash{\begin{tabular}[t]{l}$n$\end{tabular}}}}%
  \end{picture}%
\endgroup%

\caption{The generator $\sigma_i$. Note that this convention is the opposite of that of \cite{birman}. \label{fig:gen}} 
\end{figure}
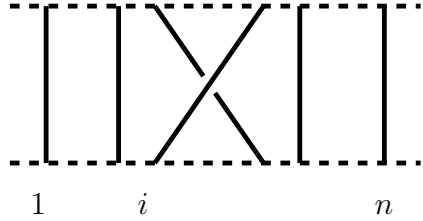

For $n = 3$, $\psi$ is defined explicitly on the generators $\sigma_1, \sigma_2$ (see Figure ~\ref{fig:gen}) as
\[ \psi(\sigma_1^{-1}) = \left[ \begin{array}{cc} -t & 1 \\ 0 & 1 \end{array}\right],  \] 
\[ \psi(\sigma_2^{-1}) = \left[ \begin{array}{cc} 1 & 0 \\ t & -t \end{array}\right].  \] 

Let $a$ be an element of $B_3$, $\hat{a}$ be the closed braid, and $e_a$ be the exponent sum of $a$. In general for $a\in B_n$, $\hat{a}$ being a knot implies that $n-1 + e_a$ is even, since if $n-1 + e_a$ is odd, a quick argument by visual inspection shows that $\hat{a}$ has more than one component. Thus for $a\in B_3$, where $\hat{a}$ is a knot, $2 \pm e_{a}$ is even and therefore $e_{a}$ is even. The Jones polynomial $J_{\hat{a}}(t)$ of $\hat{a}$ can be written in terms of $\psi$ \cite{jones}: 
\begin{align}
J_{\hat{a}}(t) &=  (-\sqrt{t})^{-e_a}(t+t^{-1} + \text{trace } \psi(a)). \label{eq:trace}
\intertext{The sign change on $e_a$ as compared to \cite[Eq. (5)]{birman}, where this form of the equation is taken from, is due to the difference in convention as indicated in Figure ~\ref{fig:gen}. When $n=3$, the Alexander polynomial of $\hat{a}$ may also be written in terms of the trace of $\psi$ \cite[Eq. (7)]{birman}: } 
(t^{-1} + 1 + t)\triangle_{\hat{a}}(t) &= (-1)^{-e_a} (t^{-e_a/2} - t^{e_a/2} \text{trace } \psi (a) + t^{e_a/2}), \label{eq:substitution}
\end{align} 
with similar adjustments on the signs.
Rearranging equations  (\ref{eq:trace}) and (\ref{eq:substitution}) above, we have the symmetric Alexander polynomial of a closed 3-braid re-written in terms of the Jones polynomial. 
\begin{equation}
(t^{-1} + 1 + t) \triangle_{\hat{a}}(t) =  (-1)^{-e_{a}} (- (-1)^{e_{a}}t^{e_{a}}J_{\hat{a}}(t) + t^{e_{a}/2+1} + t^{e_{a}/2 - 1} + t^{-e_{a} /2} + t^{e_{a} /2}). \label{eq:jtoa}  \end{equation}

This expression allows us to compute certain coefficients of the Alexander polynomial from the Jones polynomial. By Birman and Menasco's solution \cite{birman-menasco} to the classification of 3-braids, there are finitely many conjugacy classes of $B_3$, and each 3-braid is isotopic to a representative of a conjugacy class. Schreier's work \cite{schreier} puts each representative of a conjugacy class in a normal form. 

\begin{thm} \cite{schreier}  Let $b\in B_3$ be a braid on three strands, and $C$ be the 3-braid $(\sigma_1\sigma_2)^3$. Then $b$ is conjugate to a braid in exactly one of the following forms: 
\begin{enumerate}
\item \label{generic} $C^k \sigma_1^{p_1} \sigma_2^{-q_1} \cdots \sigma_1^{p_s} \sigma_2^{-q_s}$, where $k\in \mathbb{Z}$ and $p_i, q_i$, and $s$ are all positive integers, 
\item $C^k \sigma_1^p$, for $k, p\in \mathbb{Z}$, 
\item $C^k \sigma_1\sigma_2$, for $k \in \mathbb{Z}$, 
\item $C^k \sigma_1\sigma_2\sigma_1$, for $k \in \mathbb{Z}$, or
\item  $C^k \sigma_1\sigma_2\sigma_1\sigma_2$, for $k \in \mathbb{Z}$.
\end{enumerate} 
\label{thm:1}
\end{thm} 

It suffices to study the 3-braids among the conjugacy representatives above to determine which closed 3-braid is an L-space knot. 
It is straightforward to check that $C^k \sigma_1^p$ and $C^k \sigma_1\sigma_2\sigma_1$ represent links for any $k, p \in \mathbb{Z}$. Also, since $C = (\sigma_1 \sigma_2)^3$, we get that, for any $k\in \mathbb{Z}$, $C^k \sigma_1 \sigma_2$ and $C^k \sigma_1 \sigma_2 \sigma_1 \sigma_2$ represent the $(3, 3k+1)$ and $(3, 3k+2)$ torus knots, respectively. Thus we will only need to study class (\ref{generic}) of conjugacy representatives of Theorem~\ref{thm:1}. 

Recall that if a knot $K$ is an L-space knot, then the absolute value of a nonzero coefficient of the Alexander polynomial $\triangle_K(t)$ is 1, and the nonzero coefficients alternate in sign. Moreover, let $g$ be the maximum degree of $\triangle_K(t)$ in $t$, then the coefficients of the term $t^{g-1}$ is nonzero and therefore $\pm 1$. The symmetric Alexander polynomial has the two possible forms given below for an L-space knot:
\begin{align*}
t^g - t^{g-1} + \cdots + \text{terms in-between} - t^{-(g-1)} + t^{-g} \\ 
\intertext{or} 
-t^g + t^{g-1} - \cdots - \text{terms in-between} + t^{-(g-1)} - t^{-g}. 
\end{align*} 
Thus by a quick computation, we can conclude the following statement. 
\begin{lem} \label{l:productcoeff}
Suppose that $K$ is an L-space knot, then the product 
\[ \triangle_K(t) \cdot (t^{-1} + 1 + t)\]
is a symmetric polynomial with coefficients in $\{-1, 0, 1\}$, which do not necessarily alternate in sign, and the second coefficient and the second-to-last coefficient are both zero. 
\end{lem} 

The conjugacy representatives of class (1) in Theorem ~\ref{thm:1} are called \emph{generic} 3-braids. That is, $b \in B_3$ is \emph{generic} if it has the following form.
\[ b= C^k \sigma_1^{p_1} \sigma_2^{-q_1} \cdots \sigma_1^{p_s} \sigma_2^{-q_s},  \] where $p_i, q_i, k \in \mathbb{Z}$, with $p_i, q_i > 0$, and $C = (\sigma_1\sigma_2)^3$. 

The braid $a = \sigma_1^{p_1} \sigma_2^{-q_1} \cdots \sigma_1^{p_s} \sigma_2^{-q_s}$ is called an \emph{alternating 3-braid}. The first three coefficients and the last three coefficients of the Jones polynomial for this class of 3-braids, as well as the degree, are explicitly calculated in \cite{farey}. We assemble below the results we will need. 

\begin{defn} For an alternating braid $a = \sigma_1^{p_1} \sigma_2^{-q_1} \cdots \sigma_1^{p_s} \sigma_2^{-q_s}$, let 
\[ \mathbf{p} := \sum_{i=1}^s p_i, \text{ and } \mathbf{q}:= \sum_{i=1}^s q_i,\] 
so the exponent sum $e_{a} = \mathbf{p} - \mathbf{q}$. 
\end{defn} 
We have the following lemma from \cite{farey}, phrased in terms of the notations in this paper with item (c) replaced by a result within the proof. 
\begin{lem} \label{lem:fareycoeff} \cite[Lemma 6.2]{farey} Suppose that a link $\hat{a}$ is the closure of an alternating 3-braid $a$: 
\[ a = \sigma_1^{p_1} \sigma_2^{-q_1} \cdots \sigma_1^{p_s} \sigma_2^{-q_s}, \] with $p_i, q_i > 0$ and $\mathbf{p} > 1$ and $\mathbf{q} > 1$, then the following holds.
\begin{itemize}
\item[(a)] The highest and lowest powers, $M(\hat{a})$ and  $m(\hat{a})$ of $J_{\hat{a}}(t)$ in $t$ are 
\[ M(\hat{a}) = \frac{3\mathbf{q}-\mathbf{p}}{2} \text{ and } m(\hat{a}) = \frac{\mathbf{q} - 3\mathbf{p}}{2}.\] 
\item[(b)] The first two coefficients $\alpha, \beta$ in $J_{\hat{a}}(t)$ from $M(\hat{a})$ , and the last two coefficients $\beta'$, $\alpha'$ in $J_{\hat{a}}(t)$ from $m(\hat{a})$ are 
\[ \alpha = (-1)^{\mathbf{p}}, \beta = (-1)^{\mathbf{p} + 1}(s-\epsilon_{\mathbf{q}}), \beta' = (-1)^{\mathbf{q} + 1}(s-\epsilon_{\mathbf{p}}), \alpha' = (-1)^{\mathbf{q}}. \] 
where $\epsilon_{\mathbf{p}} = 1$ if $\mathbf{p} = 2$ and $0$ if $\mathbf{p} > 2$, and similarly for $\epsilon_{\mathbf{q}}$. 
\item[(c)] \cite[Eq. (14) in the proof of Lemma 6.2]{farey} Let $\gamma, \gamma'$ denote the third and the third-to-last coefficient of $J_{\hat{a}}(t)$, respectively.
 We have
\begin{align*} 
(-1)^{\mathbf{p}} \gamma &= \frac{s^2+3s}{2} - \#\{i:p_i =1\} - \# \{i:q_i=1\} - \delta_{\mathbf{q}=3}, 
\intertext{and}
(-1)^{\mathbf{q}} \gamma' &= \frac{s^2+3s}{2} - \#\{i:p_i =1\} - \# \{i:q_i=1\} - \delta_{\mathbf{p}=3},  
\end{align*}
where $\delta_{\mathbf{q}=3}$ is zero if $\mathbf{q}\not=3$ and 1 otherwise, and $\delta_{\mathbf{p}=3}$ is similarly defined.  
\end{itemize} 
\end{lem} 
Note that item (c) implies the original statement in item (c) of \cite[Lemma 6.2]{farey}.

The next result writes the Jones polynomial of a generic 3-braid in terms of the Jones polynomial of an alternating braid. Again, the statement is phrased in terms of the notations in this paper. 
\begin{lem} \label{lem:fareygen} \cite[Lemma 6.3]{farey} If $b$ is a generic braid of the form 
\[ b = C^k a,\] where $a$ is an alternating 3-braid, and let 
$J_{\hat{b}}(t)$ denote the Jones polynomial of the closure $\hat{b}$, then
\[ J_{\hat{b}}(t) = t^{-6k}J_{\hat{a}}(t) + (-\sqrt{t})^{-e_a} (t+ t^{-1})(t^{-3k} - t^{-6k}).  \]  
\end{lem} 

By equation (\ref{eq:jtoa}), we have 
\begin{align} \notag
&\triangle_{\hat{b}}(t) \cdot (t^{-1} + 1 + t) \\ \notag 
&=  (-1)^{-e_{b}} (- (-1)^{e_{b}}t^{e_{b}}J_{\hat{b}}(t) + t^{e_{b}/2+1} + t^{e_{b}/2 - 1} + t^{-e_{b} /2} + t^{e_{b} /2}). \notag 
\intertext{Putting this together with Lemma ~\ref{lem:fareygen} and noting that $e_b$ must be even in order for $\hat{b}$ to be a knot, we get the following equation.} \notag
\triangle_{\hat{b}}(t)  &=- t^{e_{b}}(t^{-6k}J_{\hat{a}}(t) + (-\sqrt{t})^{-e_a} (t+ t^{-1})(t^{-3k} - t^{-6k})) \\  \notag 
&+ t^{e_{b}/2+1} + t^{e_{b}/2 - 1} + t^{-e_{b} /2} + t^{e_{b} /2}. \notag 
\intertext{Since $b = C^ka$ and $e_b = 6k + e_a$, we have} \notag
&\triangle_{\hat{b}}(t) (t^{-1}+1+t) \\  \notag
&= - t^{(6k + e_a)}(t^{-6k}J_{\hat{a}}(t) + (-\sqrt{t})^{-e_a} (t+ t^{-1})(t^{-3k} - t^{-6k})) \\ \notag 
&+  t^{(6k + e_a)/2+1} + t^{(6k + e_a)/2 - 1} + t^{-(6k + e_a) /2} + t^{(6k + e_a) /2}  \\ 
&= - t^{e_a}J_{\hat{a}}(t) + t^{e_a/2 -1} + t^{e_a/2+1} + t^{-3k - e_a /2} + t^{3k + e_a /2}.  \label{eq:final}
\end{align}

We are now ready to determine which closed 3-braids are L-space knots. The proof is outlined as follows:  We will first rule out the cases where $s$, the length of the alternating braid $a$, satisfies $s \geq 3$. Then we rule out the case where $s = 2$.  Lemmas \ref{l:productcoeff}, \ref{lem:fareycoeff}, \ref{lem:fareygen}, and Equation \eqref{eq:final} will be used throughout. Finally, we consider the case $s=1$ and classify all the L-space knots using results from \cite{Vafaee2013}. 
\begin{proof}  

~\paragraph{\textbf{The case} $\mathbf{s \geq 3}$: }
 For a closed braid $\hat{b}$ to be an L-space knot, the right hand side of Equation \eqref{eq:final} has to have nonzero coefficients in $\pm 1$. This immediately restricts $s$ to be less than 4, since $s \geq 4$ implies that there will be two coefficients $\beta, \beta'$, the second and the second-to-last, whose absolute values are greater than or equal to 4 in $J_{\hat{a}}(t)$ by (b) of Lemma ~\ref{lem:fareycoeff}. This would result in at least one coefficient whose absolute value is greater than or equal to 2 in $\triangle_{\hat{b}}(t)(t^{-1}+1+t)$, even after possible cancellations from the terms $ t^{e_a/2 -1}, t^{e_a/2+1}, t^{-3k - e_a /2}$, and $t^{3k + e_a /2}$. Similarly, if $s=3$, then $\mathbf{p}, \mathbf{q}\geq 3$, and $|\beta|$ and $|\beta'|$ are both greater than or equal to 3, each of which would need to be cancelled out by at least two separate terms of $t^{e_a/2 -1}, t^{e_a/2+1}, t^{-3k - e_a /2}$, and $t^{3k + e_a /2}$. In addition, part (c) of Lemma ~\ref{lem:fareycoeff} gives that the absolute values of the third and third-to-last coefficients $\gamma, \gamma'$ are greater than or equal to $2$, which would also need to be canceled out in the sum of the right hand side of (\ref{eq:final}). This is impossible, so $s\not= 3$. 

~\paragraph{\textbf{The case} $\mathbf{s = 2}$: }
 If $\mathbf{p} > 2$, then from Lemma \ref{lem:fareycoeff} (b), the second-to-last coefficient $\beta' = (-1)^{\mathbf{q}+1} (s - \epsilon_{\mathbf{p}}) = (-1)^{\mathbf{q}+1} 2 $, and the third coefficient is given by
 \begin{align*} 
 \gamma &=(-1)^{\mathbf{-p}}\left(  \frac{s^2+3s}{2} - \#\{i:p_i =1\} - \# \{i:q_i=1\} - \delta_{\mathbf{q}=3} \right). 
 \intertext{If $\mathbf{p} > 2 \text{ and } \mathbf{q} = 3$, then}
 |\gamma| &\geq 5 - 1 -1 -1 = 2. \\ 
 \intertext{If $\mathbf{p} > 2 \text{ and } \mathbf{q} \not= 3$, then}
 |\gamma| & \geq 5 -1 -2 - 0 = 2. 
 \end{align*} 
  Since $\mathbf{p}, \mathbf{q}$ are either both even or both odd to make $e_a$ even, $\beta'$ and $\gamma$ have opposite signs. One of them is positive, which would not cancel out with any of the terms $ t^{e_a/2 -1}, t^{e_a/2+1}, t^{-3k - e_a /2}$, or $t^{3k + e_a /2}$. The case is similar for $\mathbf{q}> 2$ with respect to the second coefficient $\beta$  and the third-to-last coefficient $\gamma'$, so we must have that $\mathbf{p} = 2$ and $\mathbf{q}=2$. This means that $a = \sigma_1 \sigma_2^{-1} \sigma_1 \sigma_2^{-1}$. The Jones polynomial of this alternating closed braid, the figure 8 knot,  is \cite{knotinfo}
\[ J_{\hat{a}}(t) =  t^{-2}-t^{-1}+ 1-t+ t^2. \]
The figure 8 knot is known not to be an L-space knot, which one can see directly by looking at its Alexander polynomial. 
 For a braid $b = C^ka$, we have
\begin{align*}
\triangle_{\hat{b}}(t) (t^{-1}+1+t) &= -J_{\hat{a}}(t) + \frac{1}{t} + t + t^{-3k} + t^{3k} \\
&= -(t^{-2}-t^{-1}+ 1-t+ t^2) + \frac{1}{t} + t + t^{-3k} + t^{3k} \\
&=-t^{-2} + 2t^{-1} -1 + 2t - t^2 + t^{-3k} + t^{3k}.
\end{align*}
For all $k \not= 0$, this shows that the product on the left hand side of this equation has nonzero coefficients that are not $\pm 1$. This cannot happen by Lemma \ref{l:productcoeff}. This rules out the possibility that a closed 3-braid of this form can be an L-space knot. 

~\paragraph{\textbf{The case} $\mathbf{s = 1}$:}
 Assuming that both $\mathbf{p}, \mathbf{q}$ are greater than 1, 
the absolute values of the third coefficient $\gamma$ and the third-to-last coefficient $\gamma'$ of $-t^{e_{a}}J_{\hat{a}}(t)$ have the form 
\begin{align*}
& |\gamma| = \left| \left(\frac{s^2+3s}{2} - \#\{i:p_i =1\} - \# \{i:q_i=1\} - \delta_{\mathbf{q}=3}\right) \right|,  
\intertext{and}
& |\gamma'| =\left| \left(\frac{s^2+3s}{2} - \#\{i:p_i =1\} - \# \{i:q_i=1\} - \delta_{\mathbf{p}=3}\right) \right| ,
\end{align*}
 by Lemma \ref{lem:fareycoeff} (c). They are coefficients of monomials $t^{e_a}t^{(3\mathbf{q-p})/2-2} = t^{(\mathbf{q}+\mathbf{p})/2-2} $ and $t^{e_a}t^{(\mathbf{q-3p})/2+2} = t^{-(\mathbf{q}+\mathbf{p})/2 + 2}$, respectively. 
If $|\gamma|$ or $|\gamma'|$ is $\geq 2$, they need to be canceled by at least one term out of 
\[t^{e_a/2 -1}, \ t^{e_a/2+1}, \ t^{-3k - e_a /2}, \ t^{3k + e_a /2}. \] 
Plugging in $\mathbf{p-q}$ for $e_a$, we get 
\[t^{(\mathbf{p-q})/2 -1}, \ t^{(\mathbf{p-q})/2+1}, \ t^{-3k - (\mathbf{p-q}) /2}, \ t^{3k + (\mathbf{p-q}) /2}. \] 
Thus if $\mathbf{q} > 3$ and $\mathbf{p} > 1$, then $\gamma = 2$,  and $(\mathbf{q}+\mathbf{p})/2-2$ needs to equal  $(\mathbf{p}-\mathbf{q})/2-1$, $(\mathbf{p}-\mathbf{q})/2+1$, $-3k-(\mathbf{p}-\mathbf{q})/2$, or $3k + (\mathbf{p}-\mathbf{q})/2$ for the cancellation. Similarly, we have the constraints on $-(\mathbf{q}+\mathbf{p})/2 +2$ if $\mathbf{p} > 3$ and $\mathbf{q} > 1$. We examine the resulting equations and rule out values of $\mathbf{p}$ and $\mathbf{q}$ which lead to contradictions. Setting $(\mathbf{q}+\mathbf{p})/2-2$ equal to $(\mathbf{p}-\mathbf{q})/2 -1$ or $(\mathbf{p}-\mathbf{q})/2+1$ gives $\mathbf{q} = 1$ or $\mathbf{q} = 3$. Setting $(\mathbf{q}+\mathbf{p})/2-2$ equal to $-3k-(\mathbf{p}-\mathbf{q})/2$ or $3k + (\mathbf{p}-\mathbf{q})/2$ gives $\mathbf{p} = -3k+2$ or $\mathbf{q} = 3k+2$. Similarly, if $\mathbf{p} > 3$ and $\mathbf{q} > 1$, then we must have $\mathbf{p} =3$, $\mathbf{p}=1$, $\mathbf{q} = 3k+2$, or $\mathbf{p} = -3k + 2.$ 

We set aside the cases where $\mathbf{p}$ or  $\mathbf{q} \leq 3$ for now, and suppose that $k \not=0$. By the arguments made in the previous paragraph, we have that $\gamma = \gamma' = 2$, each of which needs to be canceled out by a term in  
\[ t^{-3k - (\mathbf{p-q}) /2}, \ t^{3k + (\mathbf{p-q}) /2}. \] 
Since we assume $\mathbf{p} > 3$, $\mathbf{q} > 3$, we have that $(\mathbf{q+p})/2 -2 = -3k - (\mathbf{p-q}) /2$ or $3k + (\mathbf{p-q})/2$, which means $\mathbf{p} = -3k + 2$ or $\mathbf{q} = 3k+2$. 
 We cannot have that $\mathbf{p} = -3k + 2$ \emph{and} $\mathbf{q} = 3k + 2$ since they are both supposed to be positive. Therefore we suppose that $\mathbf{p} = -3k+2 $ or $\mathbf{q} = 3k + 2$. In the first case, $k$ is negative and $-(\mathbf{p+q})/2+2 = 3k + (\mathbf{p-q})/2$. In the second case, $k$ is positive and $-(\mathbf{p+q})/2+2 = -3k-(\mathbf{p-q})/2$. This means that $\gamma, \gamma'$ cancels with the coefficients of $t^{-3k-e_a/2}$ and $t^{3k+e_a/2}$ in Equation \eqref{eq:final}. Either way, we end up having, for $k <0$, 
\begin{align*}
& \triangle_{K}(t) \cdot (t^{-1} + 1 + t)  \\ 
&= \pm t^{-\frac{-3k+2+\mathbf{q}}{2}} \mp t^{-\frac{-3k+2+\mathbf{q}}{2}+1} \pm 1 \mp \cdots \pm 1 \mp t^{\frac{-3k+2+\mathbf{q}}{2}-1} \pm t^{\frac{-3k+2+\mathbf{q}}{2}} \\
&+ t^{\frac{-3k+2-\mathbf{q}}{2}-1} + t^{\frac{-3k+2-\mathbf{q}}{2}+1},  \\
\intertext{after cancelling the third and third-to-last coefficients of $J_{\hat{a}}(t)$ with $t^{-3k-e_a/2}$ and $t^{3k+e_a/2}$. (If they do not cancel we can immediately rule out the possibility that $\hat{b}$ is an L-space knot.)
For $k>0$, we similarly have }
& \triangle_{K}(t) \cdot (t^{-1} + 1 + t)  \\ 
&= \pm t^{-\frac{3k+2+\mathbf{p}}{2}} \mp t^{-\frac{3k+2+\mathbf{p}}{2}+1} \pm 1 \mp \cdots \pm 1 \mp t^{\frac{3k+2+\mathbf{p}}{2}-1} \pm t^{\frac{3k+2+\mathbf{p}}{2}} \\
&+ t^{\frac{-(3k+2)+\mathbf{p}}{2}-1} + t^{\frac{-(3k+2)+\mathbf{p}}{2}+1}. 
\end{align*}
One of the conditions on the product $\triangle_{K}(t) \cdot (t^{-1} + 1 + t)$ is that the second and the second-to-last coefficients are equal to zero. When $k < 0$, the terms $\mp t^{-\frac{-3k+2+\mathbf{q}}{2}+1} , \mp t^{\frac{-3k+2+\mathbf{q}}{2}-1} $ are the second and the second-to-last coefficients which are not zero since we assume $\mathbf{p}, \mathbf{q} > 3$. This is impossible by Lemma \ref{l:productcoeff} if $\hat{b}$ is an L-space knot,  so this case cannot happen. A similar argument applies to rule out the case when $k > 0$. 

When both $\mathbf{p},  \mathbf{q}= 3$, we have that the alternating 3-braid $a$ takes the form $\sigma_1^3 \sigma_2^{-3}$. The Alexander polynomial of this alternating 3-braid is 
\[ \triangle_a(t) = 3+\frac{1}{t^2} - \frac{2}{t} - 2t + t^2,   \] obtained by multiplying the Alexander polynomial of the trefoil by itself, since this 3-braid is a connected sum of two (right-hand and left-hand) trefoils. It is clear from the Alexander polynomial that this knot cannot be an L-space knot due to the fact that several of its nonzero coefficients are not $\pm 1$. Now we consider a generic 3-braid $b = C^ka$ with $a = \sigma_1^{3} \sigma_2^{-3}$. Since $e_a = 0$, the highest power and the lowest power of the Jones polynomial of $\hat{a}$ are $3$ and $-3$. By 
Equation \eqref{eq:final}, 
\[ \triangle_{\hat{b}}(t)(t^{-1}+1+t) =  -J_{\hat{a}}(t) + \frac{1}{t} + t + t^{-3k} + t^{3k},\]  
where 
\[ J_{\hat{a}}(t) = 3 - \frac{1}{t^3} + \frac{1}{t^2} - \frac{1}{t} - t + t^2 - t^3.\] 
When $k \not= 0$, it is clear that the constant term 3 of $J_{\hat{a}}(t)$ will not be canceled out by the terms $\frac{1}{t}, t, t^{-3k},$ or $t^{-3k}$. Thus none of the closures of braids of the form $C^k\sigma_1^3 \sigma_2^{-3}$ will be an L-space knot. We may also rule out the case $\mathbf{p}$ or  $\mathbf{q} = 2$ since this would give a link rather than a knot. Thus, the remaining generic 3-braids whose closure can be an L-space knot are given below. 

\begin{figure}[ht]
\begin{center} 
\begin{tabular}{|c|c|}
\hline
$C^k\sigma_1^1\sigma_2^{-q}$ & for $q$ odd.  \\ 
\hline
$C^k\sigma_1^p\sigma_2^{-1}$ & for $p$ odd. \\ 
\hline
\end{tabular}
\end{center} 
\end{figure}

We now claim that $C^k\sigma_1^p\sigma_2^{-1}$, for $p $ odd and $k>0$, represents an L-space knot. Given $ (\sigma_1 \sigma_2 \sigma_1)^2 = (\sigma_2 \sigma_1)^3$ by the braid relations $\sigma_1\sigma_2\sigma_1 = \sigma_2\sigma_1\sigma_2$, note that:
\begin{align*}
	(\sigma_1 \sigma_2 \sigma_1)^{2k} \sigma_1^p \sigma_2^{-1} & \sim (\sigma_2 \sigma_1)^{3k} \sigma_1^{p} \sigma_2^{-1}\\
		&\sim  (\sigma_2 \sigma_1)^{3k-1} \sigma_1^{p+1}. 
\end{align*}

The latter braid is the twisted torus knot, $K(3, 3k-1; 2, 1)$, which is known to be an L-space knot \cite[Corollary~3.2]{Vafaee2013}. Now if $k<0$ then
\begin{align*}
	(\sigma_1 \sigma_2 \sigma_1)^{2k} \sigma_1^p \sigma_2^{-1} & \sim (\sigma_1^{-1} \sigma_2^{-1} \sigma_1^{-1}\sigma_1^{-1} \sigma_2^{-1} \sigma_1^{-1})^{-k} \sigma_1^{p}\sigma_2^{-1} \\
               &\sim \sigma_2^{-1} (\sigma_1^{-1} \sigma_2^{-1} \sigma_1^{-1}\sigma_1^{-1} \sigma_2^{-1} \sigma_1^{-1})^{-k} \sigma_1^{p}\\
               &\sim \sigma_2^{-1} (\sigma_1^{-1} \sigma_2^{-1} \sigma_1^{-1}\sigma_2^{-1} \sigma_1^{-1} \sigma_2^{-1})^{-k} \sigma_1^{p}\\
               &\sim \sigma_2^{-1} (\sigma_1^{-1} \sigma_2^{-1} \sigma_1^{-1}\sigma_2^{-1} \sigma_1^{-1} \sigma_2^{-1})^{-k} \sigma_1^{-1}\sigma_1\sigma_1^{p+1}\\
		&\sim  (\sigma_2^{-1} \sigma_1^{-1})^{-3k+1} \sigma_1^{p+1}\\
                &\sim  (\sigma_1 \sigma_2)^{3k-1} \sigma_1^{p+1}. 
\end{align*}
Using \cite[Corollary~3.2]{Vafaee2013}, we get that the closure of the latter braid represents an L-space knot. We should point out that \cite[Corollary~3.2]{Vafaee2013}, as stated, only holds for positive knots.  However, it turns out that every twisted $(p, q)$ torus knot, where the twisting happens between $p-1$ strands, admits an L-space surgery from the proof of \cite[Theorem~3.1]{Vafaee2013}. A similar argument shows that the closure of $C^k\sigma_1^1\sigma_2^{-q}$, for $q$ odd, also represents an L-space knot. Notice that in this case the knot is isotopic to the closure of $(\sigma_2 \sigma_1)^{1-3k} \sigma_1^{-p-1}$, 
so its mirror image admits a positive L-space surgery.

\end{proof} 

\bibliographystyle{amsalpha2}
\bibliography{Reference}

\end{document}